\documentclass{amsart}

\usepackage[dvipsnames]{xcolor}
\usepackage{tikz}
\usepackage{amsmath,bm,bbm,amsthm, amssymb}
\usepackage{fullpage}
\usepackage{color}
\usepackage{dsfont}
\usepackage{animate}
\usepackage{etoolbox}
\usepackage{comment}
\usepackage{pgf}
\usetikzlibrary{calc}
\usetikzlibrary{patterns}
\usetikzlibrary{arrows}
\usetikzlibrary{decorations.pathreplacing}
\usepackage[utf8]{inputenc}
\usepackage{pgfplots}

\usepackage{tcolorbox}
\usepackage{adjustbox}
\usepackage[final]{pdfpages}
\usepackage[ruled,vlined,linesnumbered,algosection,resetcount]{algorithm2e}
\usepackage{blkarray}
\usepackage{arydshln}
\usepackage{wrapfig}
\usepackage[font=scriptsize, labelfont=sc]{caption}

\theoremstyle{plain}
\newtheorem{theorem}{Theorem}
\newtheorem{proposition}[theorem]{Proposition}

\newtheorem{lemma}[theorem]{Lemma}

\theoremstyle{definition}

\newtheorem{example}[theorem]{Example}

\theoremstyle{remark}
\newtheorem{remark}[theorem]{Remark}

\def\Del{\ms{Del}_L}

\def\La{\Lambda}
\def\Ld{L}

\def\Piv{\ms{Piv}}
\def\P{\mathbb P}

\def\R{\mathbb R}
\def\EE{\mc E}
\def\QQ{\mc Q}
\def\UU{\mc U}
\def\G{\Gamma}

\def\Var{\ms{Var}}

\def\Xs{X^{\ms s}}
\def\Xps{\Xs}
\def\Inf{\ms{Inf}}
\def\Infs{\ms{Inf}^ {Y}}
\def\Infd{\ms{Inf} ^{X}}

\def\Z{\mathbb Z}

\def\been{\begin{enumerate}}
\def\bee{\begin{example}}
\def\beit{\begin{itemize}}
\def\bel{\begin{lemma}}
\def\bepr{\begin{proposition}}
\def\bep{\begin{proof}}
\def\bet{\begin{theorem}}

\def\co{\colon}
\def\da{\downarrow}
\def\de{\delta}

\def\dist{\ms{dist}}

\def\d{{\rm d}}
\def\enen{\end{enumerate}}
\def\ene{\end{example}}
\def\enit{\end{itemize}}
\def\enl{\end{lemma}}
\def\enpr{\end{proposition}}
\def\enp{\end{proof}}
\def\ent{\end{theorem}}

\def\es{\varnothing}

\def\ff{\infty}
\def\f{\frac}

\def\im{\item}
\def\lac{\la_c}

\def\coa{\ms{coarse}}
\def\fin{\ms{fine}}
\def\lad{\la_{\ms{Del}}}

\def\la{\lambda}

\def\lrsa{\leftrightsquigarrow}
\def\mc{\mathcal}
\def\ms{\mathsf}

\def\pa{\partial}

\def\sm{\setminus}

\def\tff{\uparrow\infty}

\def\th{\theta}

\def\bef{\begin{figure}[!h]}
\def\enf{\end{figure}}

\begin{document}

\title{Sharp phase transition for Cox percolation}

\author{Christian Hirsch}
\address[Christian Hirsch]{Department of Mathematics\\ Aarhus University \\  Ny Munkegade, 118, 8000, Aarhus C,  Denmark.}
\address[Christian Hirsch]{Bernoulli Institute for Mathematics, Computer Science and Artificial Intelligence \\ University of Groningen (Univ Groningen)\\ Nijenborgh 9, NL-9747 AG Groningen, Netherlands.}
\address[Christian Hirsch]{CogniGron (Groningen Cognitive Systems and Materials Center) \\ University of Groningen (Univ Groningen) \\ Nijenborgh 4, NL-9747 AG Groningen, Netherlands}
\email{hirsch@math.au.dk}
\author{Benedikt Jahnel}
\address[Benedikt Jahnel]{Weierstrass Institute for Applied Analysis and Stochastics, Mohrenstra\ss e 39, 10117 Berlin, Germany}
\email{benedikt.jahnel@wias-berlin.de}
\author{Stephen Muirhead}
\address[Stephen Muirhead]{School of Mathematics and Statistics, University of Melbourne, Melbourne, Australia}
\email{smui@unimelb.edu.au}

\begin{abstract}
We prove the sharpness of the percolation phase transition for a class of Cox percolation models, i.e., models of continuum percolation in a random environment. The key requirements are that the environment has a finite range of dependence and satisfies a local boundedness condition, however the FKG inequality need not hold. The proof combines the OSSS inequality with a coarse-graining construction.
\end{abstract}
\maketitle


\section{Introduction}
\label{int_sec}

The field of continuum percolation deals with the existence and properties of giant connected components of a geometric graph on a stochastic system of points scattered at random in Euclidean space~\cite{cPerc}. Since its early days, continuum percolation has attracted attention from researchers in wireless communication~\cite{gilbert}. This appeal is based on the prospect of using the asymptotic theory to predict the behavior of large systems of devices that interact in a peer-to-peer fashion.

\begin{wrapfigure}{r}{0.405\textwidth}
\input{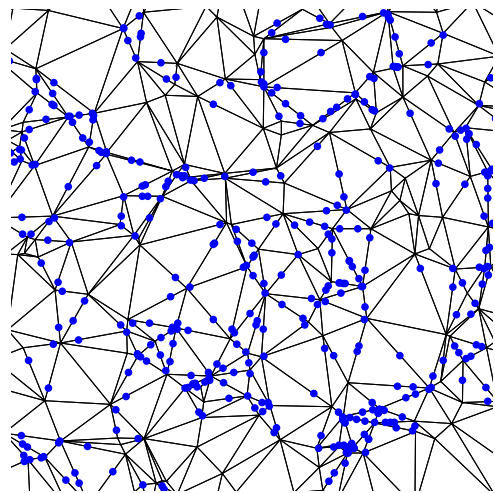}
\captionsetup{margin={0.2cm, 0.2cm}}
\caption{Devices (blue) scattered at random on the edges of a Poisson--Delaunay triangulation.}
\label{fig1}
\end{wrapfigure}

So far, the majority of results assume that the devices are scattered entirely at random in the infinite Euclidean plane, in the sense that they form a homogeneous Poisson point process. However, this assumption is in stark contrast with the topology of modern wireless networks, where devices are located predominantly on the streets of large cities, see Figure~\ref{fig1}. This discrepancy has motivated research in the direction of \emph{Cox percolation}, which can be thought of as a model of continuum percolation where the node distribution is governed by a random environment~\cite{coxPerc}. Another option is to allow the connection radii to be governed by the random environment~\cite{tykesson}.

Motivated by the need to extend the basic findings from percolation theory to random environments, \cite[Theorems 2.4, 2.6]{coxPerc} developed general conditions on the random environment ensuring a non-trivial phase transition. However, any progress beyond these basic findings was limited by the fact that the vast majority of environments relevant for applications do not satisfy the FKG inequality, which is a basic building block in percolation theory.

Recently, it has been discovered \cite{dec_tree} that the OSSS-inequality (O'Donnell, Saks, Schramm \& Servedio) is a powerful tool for analyzing percolation systems with complex spatial correlations \cite{osss2,muirhead}. Since the OSSS inequality does not rely on the FKG inequality, it is particularly attractive for analyzing Cox percolation.

Our main result shows how to apply the OSSS inequality to establish the sharpness of the phase transition in a Cox percolation model built on a random environment subject to a set of general conditions. The main features that we require from the environment are a `factor of iid' representation with finite range of dependence, and uniform local boundedness of the node intensity. We always work in the annealed model where the percolation probabilities average both over the environment and the particle placement. As a prototypical illustration of the general methodology, we apply our results to models where nodes are distributed at random on the edge set of a planar Delaunay triangulation.

As announced, the key tool to establish the main result is the OSSS inequality. However, in the context of Cox percolation, the standard approaches from the literature do not apply immediately since there are two sources of randomness: (i) the random configuration of nodes, and (ii) the random environment. Applying the OSSS inequality leads to a variance bound involving influences with respect to both sources of randomness. In contrast, when applying Russo's formula for the derivative of the percolation probability, only influences with respect to random node locations appear. In order to convert one type of influence to another, we will adapt a coarse-graining strategy from \cite{severo}.

The rest of the manuscript is organized as follows. In Section \ref{mod_sec}, we define a general framework for `factor of iid' representations of Cox processes, and provide Delaunay-based examples that are covered by this framework. Next, in Section \ref{res_sec}, we state the main results, namely Proposition \ref{pt_prop} and Theorem \ref{sharp_thm} on the non-triviality and the sharpness of the phase transition. Finally, these results are proven in Section \ref{pt_sec} and \ref{sharp_sec} respectively.

\section{Factor of iid representations of Cox processes}
\label{mod_sec}

Recall that a \textit{Cox point process} is a Poisson point process whose intensity measure (`environment') is a random Borel measure on $\R^d$. In this section, we define a general class of Cox point processes in which both, the environment and the Poisson points, may be `locally' constructed from (independent) iid processes in the background. This representation will be crucial in the proof of our main results, where we will also assume stronger properties such as finite-range dependence and uniform boundedness.

The framework is broad enough to cover environments that are absolutely continuous with respect to the Lebesgue measure, e.g., supported on a random closed set of full dimension (see Remark \ref{width_rem}), and also singular environments, e.g., supported on lower-dimensional structures such as hyperplanes or line-segments. In Section \ref{exp_sec} below, we give examples of the both types motivated by wireless communications.

\subsection{Definition}
First we introduce a class $\mathfrak{F}$ of environments which possess a `factor of iid' representation. We fix a large scale $M \ge 1$ and a fine scale $b$ such that $b^{-1}$ is an integer exceeding $2dM$; the scale $M$ will later encode the dependency range of the environment, whereas the environment will be constructed on the scale $Mb$. Then, a random Borel measure $\EE$ on $\R^d$ is in class $\mathfrak{F}$ if it is of the form
\[ \EE(\cdot) := \sum_{x \in b\Z^d} \EE_x(\cdot;Y) , \quad  \EE_x(\cdot;Y) := \UU_x(Y) \QQ_x(\cdot; Y) , \]
 where
\been
\im[(i)] $Y := \{Y_z\}_{z \in \Z^d}$ is an iid family of random elements taking values in some measurable space, 
\im[(ii)]  $\UU_x:=\UU_x(Y)$ is a random variable with values in $[0, \infty)$,
\im[(iii)] $\QQ_x(\cdot) := \QQ_x(\cdot; Y)$ is a random probability measure in the cube $Q(x; b, M):=Mx + [0, Mb]^d$, and
\im[(iv)] $\QQ_x(\cdot, Y)$ and $\UU_x(Y)$ are translation covariant, i.e., $\QQ_{x+z_0}(\{Y_{z + z_0}\}_{z \in \Z^d}) = \QQ_x(\{Y_z\}_{z \in \Z^d})$ for any $z_0 \in\Z^d$ and $x \in b\Z^d$, and similarly for $\UU_x$.
\enen
 The intuition behind this construction is that $\QQ_x$ and $\UU_x$ encode respectively the distribution of the locations and the mean number of particles in the cube $Q(x; b, M)$.

%
%
For an environment $\EE \in \mathfrak{F}$ and intensity parameter $\lambda > 0$, we next construct a Cox point process $\Xs$ with intensity $\lambda \EE$. To this end we fix a bimeasurable bijection $\G\co [0, Mb]^d \to [0, 1]$, and note that the push-forward $\G_*(\QQ_x(\cdot; Y))$ is a random probability measure on $[0, 1]$ whose cumulative distribution function will henceforth be denoted by $\Phi_x(t) := \Phi_x(t; Y)$. Next, we let 
\begin{align}
	\label{filt0_eq}
 V := \{V_x\}_{x \in b\Z^d} := \{(V_{x, i}, U_{x, i})_{i \ge 1}\}_{x \in b\Z^d}
\end{align}
be a family of iid homogeneous Poisson point processes on $[0,1] \times [0,\infty)$ with intensity $\la$. Then, we define the point process 
\begin{equation}
\label{e:xs}
\Xs := \Psi(V, Y) := \bigcup_{x \in b\Z^d}\Psi_x(V_x; Y),
\end{equation}
where the configuration of $\Xs$ inside $Q(x; b, M)$ is given as 
\begin{align}
	\label{filt_eq}
	\Psi_x(V_x; Y) := \big\{Mx + \G^{-1}(\Phi_x^{-1}(V_{x, i}; Y)) \co U_{x, i} \le \UU_x(Y) \big\}.
\end{align}
That is, $\Xs \cap Q(x; b, M)$ consists of all shifted and transformed points $\G^{-1}(\Phi_x^{-1}(V_{x, i}; Y))$ with mark $U_{x, i}$ at most $\UU_x(Y)$. Here, $\Phi_x^{-1}(v) := \inf_{y \ge 0}\{\Phi_x(y) \ge v\}$ denotes the inverse distribution function. 

One observes that, conditionally on $Y$, $\Xs \cap Q(x; b, M)$ is a Poisson point process with intensity $\la \UU_x \QQ_x$, and hence $\Xs$ is a Cox process with intensity $\la \EE$. 

\begin{remark}[Stationarity]
By construction, the Cox process $\Xs$ defined above is invariant under translations of the lattice $M\Z^d$. In certain examples (see Section \ref{exp_sec}) it may also be invariant under translations of $\R^d$.
\end{remark}

\begin{remark}[Full-dimensional environments]
\label{width_rem}
Encoding the Cox point process via the bimeasurable bijection $\Gamma$ is a bit cumbersome from a technical perspective but it allows us to cover both non-singular and singular environments simultaneously. In the former case, it may be more natural to rely on simpler alternative constructions. For instance, for environments that are uniformly distributed on a random closed subset of full dimension (as in the example in Section \ref{width_sec} below), one could  define $\EE_x(Y)$ to be the intersection of this random set with $Q(x; b, M)$, and then let $\Xs \cap Q(x; b, M)$ be given by $X_x \cap \EE_x(Y)$, where $X_x$ is a homogeneous Poisson point process in $Q(x; b, M)$ with intensity $\la > 0$.
\end{remark}

\subsection{Examples}
\label{exp_sec}
We next present examples motivated by wireless communications networks; in these examples $\Xs$ encodes the location of devices placed on a street system embedded in $\R^2$. 

\subsubsection{Delaunay network}
\label{del}
In our first example, devices are placed uniformly on the edges of a Delaunay triangulation formed from an underlying independent Poisson point process. More precisely, let $\mathcal{P}$ be a homogeneous Poisson point process on $\R^2$ with intensity $\lad > 0$, and let $\ms{Del}(\mathcal{P})$ denote the Delaunay triangulation with vertices given by $\mathcal{P}$. Then, for an intensity parameter $\la > 0$, we consider a Cox point process $\Xs$ with intensity $\lambda \EE :=  \la \textrm{Leb}(\ms{Del}(\mathcal{P})) $, where $\textrm{Leb}(\ms{Del}(\mathcal{P}))$ denotes the 1D-Lebesgue measure on~$\ms{Del}(\mathcal{P})$.

It is easy to see that $\EE \in \mathfrak{F}$, and so one can construct $\Xs$ using the general framework introduced above. More precisely, fix $M \ge 1$ and $b \le 1/4$ such that $b^{-1}$ is an integer, and for $z \in \Z^2$ define $Y_z :=  [0, Mb]^2 \cap (\mathcal{P} - Mbz)$, so that $Y:=\{Y_z\}_{z\in \Z^2}$ is an iid family of Poisson point processes and $\mathcal{P}=\bigcup_{z\in \Z^2}Y_z$. Then, $\QQ_x(\cdot; Y)$ and $\UU_x(Y)$ are respectively defined as the uniform distribution on $ \ms{Del}(\mathcal{P}) \cap Q(x; b, M) $ and the total length of $ \ms{Del}(\mathcal{P}) \cap Q(x; b, M)$. Further, $\Xs$ may be defined as in \eqref{e:xs}.

While this is a natural model for a wireless network, for our purposes it has two major drawbacks:
\begin{enumerate}
	\item[(a)] The network has infinite range of dependence.
	\item[(b)] Since we view the streets as having width $0$, the environment is singular, and there is no deterministic upper bound for the total intensity of $\Xs$ in a finite sampling window.
\end{enumerate}
This motivates us to introduce the following variants of the Delaunay network.

\subsubsection{Delaunay network superimposed with sparse grid}
\label{super}
To obtain a network with finite-range dependence, we superimpose the underlying Poisson point process with a sparse grid. More precisely, defining a large parameter $\Ld \ge 1$, and setting $M = M' L$ for a positive integer $M'$, we let $\mathcal{P}$ and $Y$ be as before, but replace $\ms{Del}(\mathcal{P})$ in the definition of $\QQ$ and $\UU$ with $\Del(\mathcal{P}) := \ms{Del}(\mathcal{P} \cup \Ld \Z^2)$. With this change, we then construct $\Xs$ in an identical manner. 

Since every triangle in the Delaunay tessellation $\Del(\mathcal{P})$ has diameter at most $cL$ for a suitable $c \ge 1$, by choosing $M'$ sufficiently large we ensure that the construction of the environment is $1$-dependent (on the scale $M = M' L$), i.e., setting $I^+(z) := z + \{-1, 0, 1\}^2$ we have for all $x\in b\Z^2$ that $\QQ_x(\{Y_{z'} \}_{z' \in \Z^2}) = \QQ_x(\{Y_{z'}' \}_{z' \in \Z^2})$ if $Y_{z''} = Y_{z''}'$ for every $z'' \in I^+(x/b)$. 

%
%
\subsubsection{Delaunay network of edges with positive width}
\label{width_sec}
To obtain a model with bounded intensity, one option is to consider the streets as having a non-zero `thickness' $w_0 > 0$. Precisely,  let $\mathcal{P}$, $Y$ and $L$ be as before, and define the random closed set $\ms{Del}^{w_0}_{L}(\mathcal{P}):= \big\{x \in \R^2\co \dist(x, \Del(\mathcal{P})) \le w_0\big\}$. Then replace $\ms{Del}(\mathcal{P})$  in the definition of $\QQ$ and $\UU$ with $\ms{Del}^{w_0}_{L}(\mathcal{P})$. By construction, the intensity of $\Xs$ is bounded by $\lambda$.

%
%
\subsubsection{Delaunay network of edges with capped density}
\label{cap_sec}
A second option to obtain a model with bounded intensity, but which retains the singular street structure, is to enforce a cap on the total intensity in any cube $Q(x; b, M)$. More precisely, let $\mathcal{P}$, $Y$ and $L$ be as before, and let $\rho > 0$ be a parameter. Then, define $\QQ_x$ to be the uniform distribution on $ \Del(\mathcal{P}) \cap Q(x; b, M) $, and define $\UU_x :=  \rho \wedge |\Del(\mathcal{P}) \cap Q(x; b, M)|$, i.e., we first measure the edge length of $\Del(\mathcal{P})$ in $Q(x; b, M)$ and then manually cap the resulting length at $\rho > 0$. By construction, the total intensity of $\Xs$ in the cube $Q(x; b, M)$ is bounded by $\la \rho$.

\begin{remark}
The Delaunay network in Section \ref{del} is stationary with respect to translations in $\R^2$, whereas after the introduction of the sparse grid $L \Z^2$ the model is only stationarity with respect to shifts in $L \Z^2$. One option to enforce the $\R^2$-stationarity could be to replace $L \Z^2$ with $L( \mathcal{V} + \Z^2)$, where $\mathcal{V}$ is uniformly distributed on $[0,1]^2$. It could also be interesting to work with a finite-range model that is intrinsically $\R^2$-stationary.
\end{remark}

%
%
\section{Main result}
\label{res_sec}
In this section, we state our main result (Theorem \ref{sharp_thm}) on the sharpness of the phase transition for Cox percolation models, i.e., the continuum percolation model built from a Cox point process $\Xs$.

First, let us be precise about the definition of the Cox percolation model. There are two equivalent ways to proceed: we can view the model as the subset of $\R^d$ formed by the union of balls of radius $1/2$ centred at the points $\Xs$; equivalently, we can consider the random geometric graph whose vertex set is $\Xs$ and whose edge-set contains all pairs of points in $\Xs$ at mutual distance less than $1$. For concreteness we will work with the former definition. We say that the model \textit{percolates} if it has an infinite connected component.

%
%
\subsection{Conditions on the environment}
\label{cond_sec}
In order to state our result, we introduce a set of conditions on the environment.

We assume that $\Xs$ is defined as in Section \ref{mod_sec} for an intensity parameter $\la > 0$ and an environment $\EE = \EE(Y) \in \mathfrak{F}$. This implies in particular that $\EE$ (and corresponingly $\Xs)$ is invariant under translations by $M\Z^d$. Moreover, we assume the dependence of $\EE$ on $Y$ is finite-range, and that $\EE$ has uniformly bounded local intensity. Precisely, define $I_b(z):= b\Z^d \cap (z + [0, 1)^d)$ and recall that $I^+(z) = z + \{-1, 0, 1\}^d$, we assume 
\been
\im[(i)] {\bf 1-dependence}, i.e., for all $z \in \Z^d$ and $x\in I_b(z)$ we have that $\EE_x(\{Y_{z'} \}_{z' \in \Z^d}) = \EE_x(\{Y_{z'}' \}_{z' \in \Z^d})$ if $Y_{z''} = Y_{z''}'$ for every $z'' \in I^+(z)$, and
\im[(ii)] {\bf uniformly bounded local intensity}, i.e., there exists a $\rho > 0$ such that $\UU_x \le \rho$ for every $x \in b\Z^d$.
\enen

%
%
Next, we impose some natural connectivity conditions on the environment $\EE$. Observe that, if the 1-neighborhood of the support of $\EE$ does not percolate, then, there is also no chance of percolation in the Cox model. To make this precise, we say that a site $x \in b\Z^d$ is \emph{non-empty} if $\UU_x(Y) > 0$. Analogously, we say that $z \in \Z^d$ is non-empty if $x$ is non-empty for some $x\in  I_b(z)$; this defines a finitely dependent site percolation model on $\Z^d$. We recall from \cite[Theorem 0.0]{domProd} that there exists $q_0(d) \in [0, 1]$ such that any 6-dependent site percolation model on $\Z^d$ with marginal probability at least $q_0(d)$ percolates, where the 6-dependence is chosen for convenience of the proofs. Hence we assume in the following
\been
\im[(iii)] {\bf coverage}, i.e., it holds that $\P\big(o \in \Z^d \text{ is non-empty}\big) > q_0(d)$.
\enen

Finally, we also need a more delicate connectivity condition that is described by the notion of \emph{essential connectedness} \cite{aldous}. To make this precise,  two sites $x, x' \in b\Z^d$ are \emph{adjacent} if they are at $d_\infty$-distance $b$. For $\eta > 0$, we say that $x \in b\Z^d$ is \emph{$\eta$-supported} if $\UU_x \ge \eta$. Introduce the enlarged cubes $I_b^+(z) := b\Z^d \cap ( z + [-1, 2)^d)$ and $I_b^{++}(z) := b\Z^d\cap( z + [-2,3)^d)$. Then we assume 
\been
\im[(iv)] {\bf essential connectedness}, i.e., that there exists $\eta > 0$ such that  with probability 1, for all $z\in \Z^d$, any non-empty $x, x'\in I_b^+(z)$ are connected by a chain of adjacent $\eta$-supported sites in $I_b^{++}(z)$.
\enen

\begin{remark}[Examples]
All the above conditions are satisfied for the Delaunay-based examples in Section \ref{width_sec} and \ref{cap_sec} by choosing $M$ of the form $M = M' L$, for $M'$ sufficiently large. On the other hand, Conditions~(i) and (ii) both fail for the example in Section \ref{del}, and Condition (ii) fails for the example in Section \ref{super}.
\end{remark}

%
\subsection{Statement of the main result}
Write $\P_\la$ for the law of $\Xs$ with parameter $\lambda$. For $A, B \subset \R^d$ we write $\{A \lrsa B  \text{ in }\Xps\}$ for the event that $A$ and $B$ are connected in the corresponding Cox percolation model, i.e., there exists a path between $A$ and $B$ in the set $\bigcup_{x \in \Xs} \{x + B(1/2)\}$, where $B(r)$ is the Euclidean ball of radius $r$. Then, defining $\Lambda_M = [-M,M]^d$, we let
$\lac := \inf\{\la >0 \co \th(\la) > 0\}$
be the critical intensity for percolation, where for $\partial \La_n := \La_{n - M} \sm \La_{n - 2M}$ we let
$$\th(\la) := \P_\la\big(\Lambda_{3M} \lrsa \ff \text{ in }\Xps\big) := \lim_{n \to \infty} \P_\la\big(\Lambda_{3M} \lrsa \partial \Lambda_n \text{ in }\Xps\big)$$
denote the probability of the percolation event.

Note that in classical continuum percolation, $\th(\la)$ is defined as $\P_\la( 0 \lrsa \ff \text{ in } \Xps)$. In the $M$-discretized setting the addition of a point at the origin is no longer natural, and hence we rely on a definition of $\th(\la)$ that is better adapted to the $M$-discretized model.

As a preliminary step, we verify that Cox percolation exhibits a non-trivial phase transition:

\bepr[Non-triviality]
\label{pt_prop}
It holds that $0 < \lac < \ff$. 
\enpr

The main result of the paper states that this phase transition is sharp, i.e., there is exponential decay of connectivity in the subcritical phase, and the percolation probability grows at least linearly in the supercritical phase.

\bet[Sharpness]
\label{sharp_thm}
The phase transition for Cox percolation is sharp:
\been
\im[(i)] $ \limsup_{n \tff}n^{-1} \log\P_\la\big(\Lambda_{3M}\lrsa \pa \Lambda_n \text{ in }\Xps\big) < 0$ holds for every $\la < \lac$, and 
\im[(ii)]   $ \liminf_{\la \da \lac }\th(\la) / (\la -\lac) > 0 $.
\enen
\ent

This generalizes a known result for the standard continuum percolation model \cite{meest,mench,zs85}, i.e., the homogeneous case in which $X^s$ is a Poisson point process.


\subsection{Possible extensions}

\subsubsection{The Delaunay network}
As mentioned, the Delaunay network model in Section \ref{del} does not satisfy the finite-range and uniform boundedness conditions, and so our result does not apply. It would be interesting to weaken these conditions so as to cover this model, but it would require new ideas.

\subsubsection{Random connectivity radii}
A natural generalization of the Cox percolation model would be to equip each point in $\Xs$ with a ball whose radius is drawn independently from a certain radius distribution, rather than a ball of radius $1/2$. If the radius distribution is bounded, the proof of Theorem \ref{sharp_thm} works unchanged, but it does not if the distribution is unbounded. In the homogeneous case where $\Xs$ is a Poisson point process, an extension to unbounded radii was achieved in \cite{raoufi_sub}, and it would be interesting to generalize this to Cox percolation.

\subsubsection{Varying connectivity radius}
For applications it may also be important to vary the radius of the connectivity instead of the intensity parameter (i.e., we fix $\lambda$ but connect all points in $\Xs$ within distance $\mu$, for varying $\mu > 0$). Unlike in the homogeneous case, these are not equivalent up to global rescaling. The proof from Theorem \ref{sharp_thm} breaks down in this case because of a lack of a Russo formula expressing the derivative of percolation probabilities as integrals over pivotal intensities (see \eqref{RusMar_eq}). However it may be possible to define an alternate notion of pivotal intensity \textit{wrt increasing the connectivity radius}, and compare this to other relevant notations of influence.

\section{Non-triviality: Proof of Proposition \ref{pt_prop}}
\label{pt_sec}

To establish the non-triviality of the percolation phase transition, we build on classical techniques to deal with finite-range dependent percolation processes, most notably the stochastic domination criterion \cite[Theorem 0.0]{domProd}.  Nevertheless, to make the manuscript more self-contained, we include some details.

%
%
\bep[Proof of Proposition \ref{pt_prop}, $\lac > 0$]
Recall the family $\{X_x\}_{x \in b\Z^d}$ of iid Poisson point processes on $[0,1] \times [0,\infty)$ with intensity $\la$ from \eqref{filt0_eq}. Recall also the uniformly bounded intensity condition, which implies that we can and will restrict each process $X_x$ to $[0,1] \times [0, \rho]$ without change to the Cox point process $X^s$. 

Call a site $z \in \Z^d$ is \emph{bad} if $X_x \ne \es$ for some $x \in  I_b(z)$, and observe that the probability that a site is bad tends to 0 as $\la \to 0$. Moreover, Cox percolation of $\Xs$ implies the percolation of the Bernoulli site percolation process of bad sites in $\Z^d$. Then the claim follows since Bernoulli site percolation has a subcritical regime \cite[Theorem 1.10]{Gr99}.
\enp

%
%
\bep[Proof of Proposition \ref{pt_prop}, $\lac < \ff$]
Call a site $x \in b\Z^d$ \emph{populated} if $\Xs \cap Q(x; b, M)  \ne \es$. Moreover, say that $z \in \Z^d$ is \emph{good} if (i) some $x \in I_b(z)$ is populated and (ii) any populated $x, x'\in I_b^+(z)$ are connected by a chain of adjacent populated sites in $I_b^{++}(z)$. First, note that percolation of good sites on $\Z^d$ implies Cox percolation of $\Xs$. Moreover, the good sites are 6-dependent, and the coverage and essential-connectedness conditions imply that the marginal probability exceeds $q_0(d)$ as $\la \to \ff$. Hence, an application of \cite[Theorem~0.0]{domProd} shows that for sufficiently large $\la > 0$, the good sites percolate.
\enp

\section{Sharpness: Proof of Theorem \ref{sharp_thm}}
\label{sharp_sec}

To prove the sharpness of the phase transition, our general strategy is to proceed in the vein of \cite{dec_tree, raoufi_sub} and rely on Russo's formula and the OSSS technique to derive a key differential inequality (see Proposition \ref{sup_crit_lem}). The particular challenges of Cox percolation are that (i) the model depends on an underlying environment, and that (ii) the FKG inequality does not necessarily hold. 

\subsection{The differential inequality}

To reflect that the model inherently depends on the scale $M$, we define 
\[\th_n(\la) := \P\big(\Lambda_{3M} \lrsa \pa \Lambda_{Mn} \text{ in }\Xps\big) , \]
setting $\th_i :=1$ for $i \le 4$ to avoid ambiguities. We will deduce the sharpness of the phase transition as a consequence of the following differential inequality, see \cite[Lemma 1.7]{raoufi_sub}.

%
%
\bepr[Differential inequality]
\label{sup_crit_lem}
Let $\la' > \lac$. Then, there exists $c_{\ms{Diff}} > 0$ such that for every $n \ge 1$ and $\la'>\la>\lac$,
$$\f{\d}{\d\la}\th_n(\la) \ge c_{\ms{Diff}} \f n{\sum_{s \le n}\th_s(\la)}\th_n(\la) (1 - \th_n(\la)).$$
\enpr

%
%
\bep[Proof of Theorem \ref{sharp_thm}]
Once Proposition \ref{sup_crit_lem} is established, we argue as in \cite{dec_tree} to show that there exists $\la_1 \in [\lac,\la']$ such that
 $\limsup_{n\tff} n^{-1}\log(\th_n(\la)) < 0$ for every $\la < \la_1$, and
 $\liminf_{\la \da \la_1}\th(\la)/ (\la - \la_1) > 0$. 
Therefore, $\lac = \la_1$.
\enp

It remains to deduce Proposition \ref{sup_crit_lem}. The key idea is to apply the OSSS inequality from \cite{osss} to the indicator of the event $\{\Lambda_{3M} \lrsa \pa \Lambda_{Mn} \text{ in }\Xps\}$. The fact that the Cox point process $\Xs$ has a `factor of iid' representation is crucial in implementing this strategy.

Recall the family $\{V_x\}_{x \in b\Z^d}$ of iid Poisson point processes from \eqref{filt0_eq}, and define the collection $X=\{X_z\}_{z\in \Z^d}$ with $X_z := \{V_x\}_{x \in  I_b(z)}$. Then, we observe that $\th_n(\la)$ may be considered as the expectation of a function $f_n(X, Y)$ on the discrete product measure $Z := \{Z_z\}_{z \in \Z^d} := \big\{(X_z, Y_z)\big\}_{z \in \Z^d}$. Applying the OSSS inequality to an algorithm $T$ determining $f_n$ gives that
\begin{equation}
	\label{osss}
	\th_n(\la) (1 - \th_n(\la)) = \Var(f_n(Z)) \le \frac{1}{2} \sum_{z \in \Z^d} \de_z(T) \Inf_z(f_n) ,
\end{equation}
where 
\been
\im[(i)] $\de_z(T) := \P(T\text{ reveals }Z_z)$ is the probability that the algorithm $T$ reveals the value of $Z_z$, and 
\im[(ii)] $\Inf_z(f_n) := \P\big(f_n(Z) \ne f_n(Z'(z)\big)$ denotes the \emph{(resampling) influence} of $Z_z$, where $Z'(z)$ is formed from $Z$ by replacing $Z_z$ with an independent copy.
\enen
By the Efron--Stein inequality,
\[ \Inf_z(f_n) \le \Infd_z(f_n) + \Infs_z(f_n), \]
where $\Infd_z(f_n) := \P\big(f_n(X,Y) \ne f_n(X'(z),Y)\big)$ and $\Infs_z(f_n) := \P\big(f_n(X,Y) \ne f_n(X,Y'(z))\big)$ denote the  analogous \emph{(resampling) influences} of $X_z$ and $Y_z$ respectively. Hence we also have 
\begin{equation}
	\label{osss_eq}
	\th_n(\la) (1 - \th_n(\la))  \le \frac{1}{2} \sum_{z \in \Z^d} \de_z(T) \big(\Infd_z(f_n) + \Infs_z(f_n)\big).
\end{equation}

%
%
Our task is to relate the right-hand side of \eqref{osss_eq} to the derivative of $\th_n(\la)$. For this we rely on Lemmas~\ref{piv_lem}--\ref{rev_lem} below, whose proofs will be given at the end of the section.

Recall that the uniformly bounded intensity condition implies that we may restrict each Poisson point process $V_x$ in \eqref{filt0_eq} to $[0,1] \times [0, \rho]$ without change to the Cox point process $\Xs$. Then, since the event $\{\Lambda_{3M} \lrsa \pa \Lambda_{Mn} \text{ in }\Xps\}$ is increasing, the infinitesimal Russo--Margulis formula gives
\begin{align}
	\label{RusMar_eq} 
\f{\d}{\d\la}\th_n(\la) = \la \sum_{x \in b\Z^d}\int_{[0, 1] \times [0, \rho]}\Piv_x(r, u) \d(r, u),
\end{align}
with
$$\Piv_x(r, u): = \P\big(f_n(V_{x,r,u}, Y) \ne f_n(X, Y)\big),$$ 
where $V_{x, r, u}$ is the collection $\{V'_{x'}\}_{x'}$ determined by $V'_{x'} = V_{x'}$ if $x' \ne x$ and $V'_x = V_x \cup \{(r, u)\}$. In words, $\Piv_x(r, u)$ is the probability that the event $\{\Lambda_{3M} \lrsa \pa \Lambda_{Mn} \text{ in }\Xps\}$ does not occur, but does occur after adding $(r, u)$ to $V_x$.

Our first lemma states that one can control the integrated infinitesimal pivotal probabilities in \eqref{RusMar_eq} in terms of the corresponding discrete influences $\Infd_z$:
%
%
\bel[Pivotality and influence]
\label{piv_lem}
Let $\la' > 0$. Then, there exists $c_{\Piv} = c_{\Piv}(\la') > 0$ such that for every $n \ge 1$, $\la \le \la'$ and $z \in \Z^d$,
$$  \Infd_z \le c_{\Piv}\la\sum_{x \in  I_b(z)}\int_{[0, 1] \times [0, \rho]}\Piv_x(r, u) \d(r, u) .$$
\enl

Our second lemma states that we can bound $\Infs_z$ in terms of $\Infd_{z'}$  for $z' \in I^{++}(z) := z + \{0, \pm 1, \pm 2\}^d$.
%
%
\bel[Poisson- and environment influences]
\label{inf_lem}
There exists $c_{\Inf} > 1$ such that for every $n \ge 1$, $\la > 0$ and $z \in \Z^d$,
$$\Infs_z \le  c_{\Inf} \sum_{z' \in I^{++}(z)}\Infd_{z'}.$$
\enl

Our final lemma bounds the revealment probabilities of a suitable randomized exploration algorithm:

%
%
\bel[Revealment probabilities]
\label{rev_lem}
For every $n \ge 16$ there exists a randomized exploration algorithm $T$ determining $f_n$ such that for every $\la > 0$ and $z\in \Z^d$,
$$ \de_z(T) \le \f8n\sum_{s \le n} \th_s(\la).$$
\enl

%
%
Let us use Lemmas \ref{piv_lem}--\ref{rev_lem} to conclude the proof of Proposition \ref{sup_crit_lem}:

\bep[Proof of Proposition \ref{sup_crit_lem}]
First, Lemmas \ref{inf_lem} and \ref{rev_lem} show that 
$$\f n{\sum_{s \le n}\th_s(\la)}\sum_{z\co  \Lambda_n \cap \Z^d}\de_z(T) (\Infd_z + \Infs_z) \le 8\sum_{z\in \Lambda_n\cap \Z^d}(\Infd_z + \Infs_z) \le8 (5^d + 1)c_{\Inf}\sum_{z\in  \Lambda_{n+2} \cap \Z^d}\Infd_z.$$
Hence, invoking Lemma \ref{piv_lem} together with the Russo--Margulis formula~\eqref{RusMar_eq} concludes the proof.
\enp

\subsection{Proof of the auxiliary lemmas}

%
%
First, we establish Lemma \ref{inf_lem}, i.e., the domination of $\Infs_z$ by a multiple of $\Infd_z$, for which we adapt a coarse-graining strategy from \cite{severo}.

\bep[Proof of Lemma \ref{inf_lem}]
We present a detailed proof in the case where $3 \le |z|_\ff \le n - 3$, noting that the arguments in the remaining cases are very similar. Let $Y'(z) = \{Y'_w\}_{w \in \Z^d}$ be the $z$-resampling of $Y = \{Y_w\}_{w \in \Z^d}$, i.e., the component of $\{Y_w\}_{w \in \Z^d}$ with index $z$ is replaced by an independent copy. Note that, by the assumption of $1$-dependence, the resampling does not modify the environment $\EE(Y)$ outside $I^+(z)$. 
We introduce a `coarse grained' version of the event $\big\{f_n(X, Y) \ne f_n(X, Y'(z))\big\}$, which depends on $X$ only through the configuration of $V_x$ for $x \not \in I_b^{++}(z)$. 
 First, we let $E_{\coa, -}$ denote the event that $\Lambda_{3M} \not\lrsa \pa \Lambda_{Mn}$ holds if $V_x   \cap \big([0, 1] \times [0, \rho]\big) = \es$ for all $x \in I_b^{++}(z)$. In particular, $E_{\coa, -}$ occurs under the event $\big\{f_n(X, Y) \ne f_n(X, Y'(z))\big\}$.
Next, we let $E_{\coa, +}$ denote the event that if $V_x \cap \big([0, 1] \times [0, \eta]\big) \ne \es$ for every  $x \in  I_b^{++}(z)$, then $\Lambda_{3M}\lrsa \pa \Lambda_{Mn}$. Now, setting $Q' :=(Mz + [-M, 2M)^d)$, we note that, under the event $\big\{f_n(X, Y) \ne f_n(X, Y'(z))\big\}$, there are particles $X_i,X_j \in \Xps \sm Q'$ such that  (i) $X_i$ connects to $\Lambda_{3M}$ outside $Q'$, and (ii) $X_j$ connects to $\pa \Lambda_{Mn}$ outside $Q'$. Hence, we conclude from the essential-connectedness condition that also the event $E_{\coa, +}$ occurs under the event $\big\{f_n(X, Y) \ne f_n(X, Y'(z))\big\}$. 

Defining $E_{\coa} := E_{\coa, -} \cap E_{\coa, +}$, we have shown that $\Infs_z \le \P(E_{\coa})$, and so it suffices to show that $ \P(E_{\coa})\le c_{{\Inf}}\sum_{z' \in  I^{++}(z)}\Infd_{z'}$. To that end, we let $E_{\fin, -} := \{V_x \cap ([0,1] \times [0,\rho])= \es \text{ for all }x \in  I_b^{++}(z)\}$. Similarly, we let $E_{\fin, +}$ denote the event that $V_x \cap ([0, 1] \times [0, \eta]) \ne \es$ for every $x \in I_b^{++}(z)$. Finally, we write $X^*$ for the Poisson point process obtained by resampling $\{V_x\}_{x \in I_b^{++}(z)}$. Then, by the independence property of Poisson point processes,
\begin{align*}
	&\P\big(f_n(X, Y) \ne f_n(X^*, Y)\big) \\
	&\quad \ge \P\big((\{V_x\}_{x \not \in I_b^{++}(z)}, Y) \in E_{\coa}, \{V_x\}_{x \in   I_b^{++}(z)} \in E_{\fin, -},\{V^*_x\}_{x \in I_b^{++}(z)} \in E_{\fin, +} \big)\\
	&\quad=  \P\big((\{V_x\}_{x \not \in I_b^{++}(z)}, Y) \in E_{\coa}\big)\P\big( \{V_x\}_{x \in   I_b^{++}(z)} \in E_{\fin, -}\big)\P\big(\{V^*_x\}_{x \in I_b^{++}(z)} \in E_{\fin, +} \big),
\end{align*}
where  the second and third factor in this product are bounded away from 0.
Finally, by the Efron--Stein inequality, we see that 
$\P\big(f_n(X, Y) \ne f_n(X^*, Y)\big) \le \sum_{z' \in  I^{++}(z)}\Infd_{z'},$
thereby concluding the proof.
\enp

Next, we prove Lemma \ref{piv_lem} through a short computation using tail estimates for Poisson random variables: 

%
%
\bep[Proof of Lemma \ref{piv_lem}]
We condition on $Y$ and note that by the superposition theorem (\cite[Theorem 3.3]{poisBook}) we may think of the collection $X_z = \{V_x\}_{x \in I_b(z)}$ as a homogeneous Poisson point process on $[0, 1] \times [0, \rho] \times I_b(z)$ with intensity $\la$. Hence, $X_z = \{P_i\}_{i \le N}$ where $N$ is a Poisson random variable with parameter $\la_* := \la \rho b^{-d}$ and the $\{P_i\}_{i \ge 1}$ are iid uniform on $[0, 1] \times [0,\rho] \times I_b(z)$. Moreover, writing $X':= \{X_{z'}\}_{z' \ne z}$, we let 
$$K:= \sup\big\{k \ge 0\co \Lambda_{3M} \not\lrsa \partial \Lambda_{Mn}\text{ in }\Psi((\{P_i\}_{i \le k} \cup X') ) \big\}$$
denote the maximum number of Poisson points that can be added to $I_b(z)$ such that the percolation event does not happen. Note that $K$ may also take the values $-\ff$ or $\ff$. The introduction of the quantity $K$ has the advantage that the pivotal probabilities can be concisely represented via 
$\P(K = N)= b^d \sum_{x \in I_b(z)} \int_{[0, 1]^2}\Piv_x(r, u) \d(r, u). $
Similarly, we can bound the influences through 
$ \Infd_z \le 2 \P(N > K \ge 0).$ Using the tail probabilities of a Poisson random variable shows that for every $k \ge 0$, 
$$\f{\P(N > k)}{\la\P(N = k)} = \sum_{\ell \ge k + 1} \f{k!}{\ell!}\la_*^{\ell - k - 1}.$$
Noting that the right-hand side is bounded by $\exp(\la_*)$ concludes the proof.
\enp

\begin{remark}
Another approach to proving Lemma \ref{piv_lem} which may be less sensitive to the Poisson assumption would be to first discretize space and then to invoke the Efron--Stein inequality to aggregate the coordinates.
\end{remark}

%

%
To finish we prove Lemma \ref{rev_lem}, i.e.,\ we describe the randomized algorithm leading to the asserted bound on the revealment probabilities. In essence, the proof can be adapted from previous results in the literature, e.g.,~\cite[Lemma 3.3]{raoufi_sub}. Nevertheless, to make the presentation self-contained, we provide a brief overview.
 
\bep[Proof of Lemma \ref{rev_lem}]
For every $6 \le m \le n - 3$ we construct an algorithm $T^m$ determining $\{\Lambda_{3M} \lrsa \partial \Lambda_{Mn}\}$ as follows. During the algorithm, a site $z\in\Z^d$ is called \emph{active} if it is revealed but the neighborhood $I^{++}(z)$ is not yet entirely revealed.
\been
\im[(i)] First, reveal $Z_z$ for all $z$ with $\big||z|_\ff - m\big| \le 3$. This determines the point configuration for all $z$ with $\big||z|_\ff - m\big| \le 2$. Let $S$ denote the union of all connected components intersecting $\partial \La_{Mm}$. 
\im[(ii)] Pick an active $z$ with $S\cap Q(x; b, M) \ne \es$ for some $x \in Q(z,b)$, reveal $I^{++}(z)$, and grow the components from $S$ with the particles from $\Xps \cap Q(x'; b, M)$ for $x' \in  I_b^+(z)$.
\im[(iii)] Continue this exploration until there is no more active $z$ with $S \cap (Mz + [0, M)^d) \ne \es$, or a connection from $\Lambda_{3M}$ to $\partial \Lambda_{Mn}$ is found.
\enen
Next, we note that the revealment probabilities $\de_z(T^m)$ are bounded above by the percolation probabilities in the sense that
$$\P(T^m \text{ reveals }z) \le \P\big((Mz + \Lambda_{3M}) \lrsa \partial \La_{Mm}\big) \le \th_{|m - |z|_\ff|}.$$
Thus, picking $m \in \{6, \dots, n - 3\}$ at random shows the asserted
$\de_z(T) \le \f2{n - 8}\sum_{ m \le n}\th_m \le \f4n\sum_{m \le n}\th_m.$
\enp

\subsection*{Acknowledgment.}
 The authors are indebted to R.~Lachi\`eze-Rey and E.~Cali for illuminating discussions and remarks on early versions of the manuscript. CH would like to acknowledge the financial support of the CogniGron research center and the Ubbo Emmius Funds (Univ.~of Groningen). BJ was supported by the German Research Foundation under Germany's Excellence Strategy MATH+: The Berlin Mathematics Research Center, EXC-2046/1 project ID: 390685689, and the Leibniz Association within the Leibniz Junior Research Group on Probabilistic Methods for Dynamic Communication Networks as part of the Leibniz Competition. SM was supported by the Australian Research Council (ARC) Discovery Early Career Researcher Award DE200101467.

\bibliographystyle{abbrv}
\bibliography{lit}

\end{document}